\documentclass{amsart}
\font\bit=cmbxti10
\font\sbit=cmbxti10 scaled 700

\def\bb{{\hbox{\bit b}}}
\def\sbb{{\hbox{\sbit b}}}

\newsymbol\diagdown 231F
\newsymbol\digamma 207A

\let\wtil=\widetilde
\let\what=\widehat
\let\\=\backslash
\let\sse=\subseteq
\let\limply=\Longrightarrow

\def\newmatrix#1{\null\,\vcenter{
		\baselineskip=8pt\mathsurround=-0pt\ialign{
		\hfil ${##}$
		\hfil &&
		\hfil ${##}$
		\hfil \crcr
		\mathstrut \crcr
		\noalign{\kern-\baselineskip}#1 \crcr
		\mathstrut \crcr
		\noalign{\kern-\baselineskip} \crcr }}\!}

\def\FF{{\mathbb F\kern.5pt}}
\def\RR{{\mathbb R\kern.5pt}}
\def\B{{\mathcal B}}
\def\J{{\mathcal J}}
\def\Le{{\mathcal L}}
\def\M{{\mathcal M}}
\def\N{{\mathcal N}}
\def\R{{\mathcal R}}
\def\Se{{\mathcal S}}
\def\X{{\mathcal X}}
\def\Y{{\mathcal Y}}
\def\Z{{\mathcal Z}}

\def\noi{\noindent}
\def\0{{\{0\}}}
\def\smallfrac#1#2{{\textstyle{\frac{#1}{#2}}}}
\def\QED{{\hfill\hfill\qed}}
\def\x{{\times}}
\def\ve{{_\vee\kern-1pt}}
\def\we{{_\wedge\kern-1pt}}
\def\hotimes{{\kern1pt\what\otimes}}
\def\otimesv{{\kern1pt\otimes_\ve}}
\def\otimesw{{\kern1pt\otimes_{\kern-1pt\we}}}
\def\hotimesv{{\kern1pt\what\otimes_\ve}}
\def\hotimesw{{\kern1pt\what\otimes_{\kern-1pt_\we}}}

\def\emap{\hbox to25pt{\rightarrowfill}}

\def\nmap{\Big\uparrow}
\def\diagdownBOX{\hbox{$\diagdown$}}
\def\searrowBOX{\hbox{\hglue6.5pt$\searrow$}}
\def\semap{\vbox{\offinterlineskip
\diagdownBOX\vglue-1pt\searrowBOX\vglue-6pt}}

\def\Remark #1{\vskip6pt\noi{\bf{Remark} #1$.$}}
\def\Proposition #1{\vskip6pt\noi{\bf{Proposition} #1$.$}}
\def\Theorem #1{\vskip6pt\noi{\bf{Theorem} #1$.$}}
\def\Corollary #1{\vskip6pt\noi{\bf{Corollary} #1$.$}}

\begin{document}

\vglue-30pt\noi
\hfill{\it Mathematica Slovaca}\/,
{\bf 72} (2022) 959--968

\vglue37pt
\title
[On Extensions of Bilinear Maps]
{On Extensions of Bilinear Maps}
\author{C.S. Kubrusly}
\address{Catholic University of Rio de Janeiro, Brazil}
\email{carlos@ele.puc-rio.br}
\renewcommand{\keywordsname}{Keywords}
\date{April 24, 2021}

\begin{abstract}
The paper deals with extension of bounded bilinear maps$.$ It gives a
necessary and sufficient condition for extending a bounded bilinear map on
the Cartesian product of subspaces of Banach spaces$.$ This leads to a full
characterization for extension of bounded bilinear maps on the Cartesian
product of arbitrary subspaces of Hilbert spaces$.$ Applications concerning
projective tensor products are also investigated.
\end{abstract}

\maketitle

\section{Introduction}

The purpose of this paper is to prove an extension result for bilinear maps$.$
This will be presented in Theorem 5.2$.$ It gives a necessary and sufficient
condition for a bounded bilinear map to be extended from the Cartesian product
${\M\x\N}$ of linear manifolds $\M$ and $\N$ of Banach spaces $\X$ and $\Y$ to
the Cartesian product ${\X\x\Y}$ of the Banach spaces$.$ Such conditions are
imposed on the linear manifolds only$.$ This leads to a complete unconditional
statement for the extension of bounded bilinear maps acting on the Cartesian
product of subspaces of Hilbert spaces$.$ Applications related to extensions
of bounded linear transformations on projective tensor prod\-ucts are
considered as well.

\vskip6pt
The paper is organized as follows$.$ Notation and terminology are set in
Section 2$.$ A brief review on the bilinear extension problem is considered
in Section 3$.$ Auxiliary results required in the sequel are brought together
in Section 4$.$ The main results are treated in Section 5, followed by an
application in Section 6.

\section{Notation and Terminology}

In the present context the terms {\it forms}\/ and {\it functionals}\/,
{\it bounded linear}\/ and {\it continuous linear}, {\it bounded bilinear}\/
and {\it continuous bilinear}, are pairwise synonyms and we use both forms
freely; and $\FF$ denotes either the real or the complex field.

\vskip6pt
All linear spaces are over the same field $\FF.$ Let ${\X,\Y,\Z}$ be nonzero
linear spaces$.$ A {\it bilinear map}\/ ${\phi\!:\X\x\Y\to\Z}$ is a function
from the Cartesian product ${\X\x\Y}$ of linear spaces to a linear space $\Z$
whose sections are linear transformations$.$ In other words, let
$\phi^y\!=\phi(\cdot,y)={\phi|_{\X\x\{y\}}\!:\X\!\to\Z}$ be the $y$-section of
the bilinear map $\phi$ and let
$\phi_x\!=\phi(x,\cdot)={\phi|_{\{x\}\x\Y}\!:\Y\to\Z}$ be the $x$-section of
$\phi.$ These functions $\phi^y$ and $\phi_x$ are linear transformations$.$
{From} now on suppose ${\X,\Y,\Z}$ are normed spaces$.$ If two normed spaces
$\X$ and $\Y$ are isometrically isomorphic, and if ${y\in\Y}$ is the
isometrically isomorphic image of ${x\in\X}$, then write ${\X\cong\Y}$ and
${x\cong y}.$ By a subspace $\M$ of a normed space $\X$ we mean a
{\it closed}\/ linear manifold of $\X$ (equipped with the norm inherited from
$\X).$ If $\M$ is a linear manifold of $\X$, then $\M^-\!$ will denote its
closure in $\X.$ A bilinear map is bounded if
$\sup_{0\ne x\in\X,\,0\ne y\in\Y}\smallfrac{\|\phi(x,y)\|}{\|x\|\,\|y\|}$ is
finite$.$ In this case set
$\|\phi\|
=\sup_{0\ne x\in\X,\,0\ne y\in\Y}\smallfrac{\|\phi(x,y)\|}{\|x\|\,\|y\|}$
so that ${\|\phi(x,y)\|\le\|\phi\|\kern1pt\|x\|\kern1pt\|y\|}$ for
\hbox{every} ${(x,y)\in\X\x\Y}.$ This defines a norm on the linear space of
all bounded bilinear maps$.$

\vskip6pt
A bilinear map is continuous (regarding the product topology in ${\X\x\Y}$) if
and only if it is bounded$.$ Let ${\bb[\X\x\Y,\Z]}$ denote the normed space of
all bonded bilinear maps ${\phi\!:\X\x\Y\to\Z}$, and let $\B[\X,\Y]$ denote
the normed space of all bounded linear transformations
${T\!:\X\kern-1pt\to\kern-1pt\Y}.$ The range of ${T\in\B[\X,\Y]}$ (notation:
$\R(T)=T(\X)=$ $\hbox{range}\,(T)\kern.5pt$) is a linear manifold of $\Y$, and
its kernel (notation: $\N(T)=T^{-1}(\0)=$ $\hbox{kernel}\,(T)\kern.5pt$) is a
subspace of $\X.$ Let ${\X^*\!=\B[\X,\FF]}$ be the dual of $\X.$ If one of
$\X$ or $\Y$ is a Banach space, then $\phi$ lies in ${\bb[\X\x\Y,\Z]}$ if and
only if $\phi^y$ lies in ${\B[\X,\Z]}$ and $\phi_x$ lies in ${\B[\Y,\Z]}.$
Both ${\bb[\X\x\Y,\Z]}$ and $\B[\X,\Z]$ are \hbox{Banach} spaces if and only
if $\Z$ is$.$ (For properties on bounded bilinear maps see, e.g.,
\cite[Section 1.1, p.8]{DF}.)

\vskip6pt
The algebraic tensor product of linear spaces $\X$ and $\Y$ is a linear
space ${\X\otimes\Y}$ associated with a bilinear map
${\theta\!:\X\x\Y\to\X\otimes\Y}$ whose range spans ${\X\otimes\Y}$ with
the following property: for every linear map ${\phi\!:\X\x\Y\to\Z}$ into
any linear space $\Z$ there exists a (unique) linear transformation
${\Phi\!:\X\otimes\Y\to\Z}$ for which the diagram
$$
\newmatrix{
\X\x\Y & \kern2pt\buildrel\phi\over\emap & \kern-1pt\Z                   \cr
       &                                 &                               \cr
       & \kern-3pt_\theta\kern-3pt\semap & \kern4pt\nmap\scriptstyle\Phi \cr
       &                                 & \phantom{;}                   \cr
       &                                 & \kern-2pt\X\otimes\Y          \cr}
$$
commutes$.$ Set ${x\otimes y=\theta(x,y)}$ for each ${(x,y)\in\X\x\Y}.$ These
are the {\it single tensors}\/$.$ An element $\digamma$ in ${\X\otimes\Y}$ is
a finite sum ${\sum_ix_i\otimes y_i}$ of single tensors$.$ (For an
expo\-sition on algebraic tensor products see, e.g., \cite{Kub2}.)

\vskip6pt
Let $\X$ and $\Y$ be normed spaces$.$ Two reasonable crossnorms on
${\X\otimes\Y}$ are the {\it injective}\/ ${\|\cdot\|_\ve}$ and
{\it projective} ${\|\cdot\|_\we}$ norms, which are given by
$$ 
\|\digamma\|_\ve
=\sup_{\|f\|\le1,\,\|g\|\le1,\;f\in\X^*\!,\,g\in\Y^*}
\Big|{\sum}_if(x_i)\,g(y_i)\Big|,
$$
 \vskip-4pt\noi
$$
\|\digamma\|_\we
=\inf_{\{x_i\},\,\{y_i\},\;\digamma=\sum_ix_i\otimes y_i}
{\sum}_i\|x_i\|\,\|y_i\|,
$$
for every ${\digamma={\sum}_ix_i\otimes y_i\in\X\otimes\Y}.$ Let
${\X\otimesv\Y}$ and ${\X\otimesw\Y}$ denote a tensor product space
${\X\otimes\Y}$ equipped with ${\|\cdot\|_\ve}$ or ${\|\cdot\|_\we}.$ Their
completions, denoted by ${\X\hotimesv\Y}$ and ${\X\hotimesw\Y}$, are referred
to as the {\it injective and projective tensor products}\/$.$ (For expositions
on the Banach spaces ${\X\hotimesv\Y}$ and ${\X\hotimesw\Y}$ see, e.g.,
\cite{Jar}, \cite{DF}, \cite{Rya}, \cite{DFS}.)

\section{Preliminaries}

A very brief review of previous results on bilinear extension under
restrictions.

\vskip6pt\noi
Let ${\X,\Y,\Z}$ be normed spaces over the same field and let $\M$ and $\N$
be linear manifolds of $\X$ and $\Y$, respectively$.$ Although there is no
Hahn--Banach Theorem for bilinear functionals, variants of it may hold if
extra assumptions are imposed$.$ Consider the general case of a bilinear map
${\phi\!:\!\M\x\N\!\to\Z}$ taking values in a Banach space $\Z$ rather than in
$\FF.$ Having in mind the extension of $\phi$ to ${\what\phi\!:\X\x\Y\to\Z}$
acting on ${\X\x\Y}$, extra assumptions may be placed on any of the four
players involved, namely, the normed spaces $\X$ and $\Y$, the linear
manifolds $\M$ and $\N$, the bounded bilinear map $\phi$, and the Banach space 
$\Z.$ The first goal in the present paper is to generalize an extension for
bounded bilinear forms proposed in \cite{Hay}, which will be discussed in
Section 5$.$ Other ways to face the bounded bilinear extension problem was to
assume ${\X=\Y}$ and ${\Z=\FF}$ (i.e., bilinear or multilinear forms rather
than maps)$.$ Quite often it is assumed ${\M=\N}.$ Along this line, the role
played by the embedding of $\M$ into $\X$ was analyzed in \cite{CGJ1}; the
case of ${\X=L^1}$ was presented in \cite{CGJ2}, \cite{DGMP}; properties of
the bilinear (multilinear) form itself were investigated in \cite{JPPV}; the
restriction operator taking bounded bilinear forms on ${\X\x\X}$ into bounded
bilinear forms into ${\M\x\M}$ was considered in \cite{CGDPS}; a
characterization of the sequence space $c_0$ in terms of extendible bilinear
forms was given \cite{CS}; and a necessary and sufficient condition for
bilinear forms to be extendible was presented in \cite{Bow} making a
connection with an integral representation of bounded bilinear functionals
identified with bounded linear functionals on the injective tensor product
(see, e.g., \cite[Proposition 1.1.21]{DFS})$.$ All these were done for
bilinear (multilinear) forms (rather than for $\Z$-valued bilinear maps),
having in mind the quest for conditions (or particularizations) ensuring
Hahn-Banach type results for continuous bilinear functionals.

\section{Auxiliary Results}

A subspace of a normed space is {\it complemented}\/ if it has a subspace as
an algebraic complement$.$ A normed space is {\it complemented}\/ if every
subspace of it is complemented$.$ If a Banach space is complemented, then it
is isomorphic (i.e., topologically isomorphic) to a Hilbert space, and so
complemented Banach spaces are identified with Hilbert spaces \cite{LT} (see
also \cite{Kal})$.$ We will need the following well-known results, one on
complemented subspaces and the other on bilinear maps.

\Proposition{4.1}
{\it Let\/ $\M$ be a subspace of an normed space\/ $\X.$ If\/ there exists a
continuous projection\/ ${E\!:\X\kern-1pt\to\X}$ with\/ ${\R(E)=\M}$, then\/
$\M$ is complemented$.$ The converse holds if\/ $\X$ is a Banach space}\/.

\proof
See, e.g., \cite[Remark A.4]{ST2} among many others. \QED

\Proposition{4.2}
{\it For an arbitrary triple\/ ${(\X,\Y,\Z)}$ of normed spaces}\/,
$$
\bb[\X\x\Y,\Z]\cong\B[\X,\B[\Y,\Z]\kern1pt].
$$

\proof
For ${\Z=\FF}$ and with $\cong$ standing for isometric isomorphism we get
$$
\bb[\X\x\Y,\FF]\cong\B[\X,\Y^*].
$$
See \cite[Section 1.4, p.9]{DF}$.$ An extension from $\FF$ to $\Z$ follows the
same argument. \QED

\vskip6pt
If $\M$ and $\N$ are subspaces of Banach spaces $\X$ and $\Y$, then
${\M\otimesv\N}$ is a linear manifold of ${\X\otimesv\Y}$, and consequently
${\M\hotimesv\N}$ is a subspace of the injective tensor product
${\X\hotimesv\Y}.$ In general, this fails for the projective norm
${\|\cdot\|_\we}$ of the projective tensor product ${\X\hotimesw\Y}.$
Indeed, as is readily verified,
$\|\digamma\|_{\X\otimesw\Y}\le\|\digamma\|_{\M\otimesw\N}$ for every
$\digamma\!$ in ${\M\otimesw\N}$, the inequality may be strict, and
${\M\otimesw\N}$ is a subspace of ${\X\otimesw\Y}$ if and only if
${\|\digamma\|_{\X\otimesw\Y}=\|\digamma\|_{\M\otimesw\N}}$ for every
$\digamma\!$ in ${\M\otimesw\N}.$ We will need the following well-known
results on projective tensor products.

\Proposition{4.3}
{\it Suppose\/ $\M$ and\/ $\N$ be complemented subspaces of Banach spaces\/
$\X$ and\/ $\Y$, respectively$.$ If ${\M\otimesw\N}$ is a linear manifold
of the normed space ${\X\otimesw\Y}$, then it is a subspace$.$ If\/
${\M=\R(E)}$ and\/ ${\N=\R(P)}$ for projections\/ $E$ in\/ $\B[\X,\X]$ and\/
$P$\/ in\/ $\B[\Y,\Y]$ with\/ ${\|E\|=\|P\|=1}$, then\/ ${\M\otimesw\N}$ is a
linear manifold of}\/ ${\X\otimesw\Y}$.

\proof
See, e.g., \cite[Proposition 2.4]{Rya}. \QED

\Proposition{4.4}
{\it For an arbitrary triple\/ ${(\X,\Y,\Z)}$ of Banach spaces}\/,
$$
\bb[\X\x\Y,\Z]\cong\B[\X\hotimesw\Y,\Z].
$$

\proof
This is the {\it universal mapping principle}\/$.$ See, e.g.,
\cite[Theorem 1.1.8]{DFS}. \QED

\Proposition{4.5}
{\it An alternate expression for the projective norm}\/:
$$
\|\digamma\|_\we
=\sup_{\|\phi\|\le1,\,\phi\in\sbb[\X\x\Y\!,\,\FF]}
\Big|{\sum}_i\phi(x_i,y_i)\Big|
\quad\;\hbox{for every}\;\quad
\digamma={\sum}_ix_i\otimes y_i\in\X\otimesw\Y.
$$

\proof
See, e.g., \cite[p.23]{Rya}$.$ This follows in part by Proposition 4.4 and
the fact that, if $\what\X$ is a completion of $\X$ and $\Z$ is a Banach
space, then ${\B[\what\X,\Z]\cong\B[\X,\Z]}.$ Hence
$$
\bb[\X\x\Y,\FF]\cong(\X\hotimesw\Y)^*\cong(\X\otimesw\Y)^*.      \eqno{\qed}
$$

\section{An Extension for Bounded Bilinear Maps}

If one imposes appropriate restrictions on the linear manifolds $\M$ and $\N$
of $\X$ and $\Y$, then we get (i) the next extension result for bounded
bilinear forms due to \cite{Hay}, and (ii) an extension of it for bilinear
maps as in the subsequent theorem.

\Theorem{5.1}
\cite[Corollary 2]{Hay}$.$
{\it If\/ $\M$ and\/ $\N$ are subspaces of\/ Banach spaces\/ $\X$ and\/ $\Y$,
respectively, and if there exists a projection of norm one of\/ $\X$ onto\/
$\M$ and a projection of norm one of\/ $\Y$ onto\/ $\N$, then every bounded
bilinear functional on\/ ${\M\x\N}$ can be extended to\/ ${\X\x\Y}$ with the
same norm}\/.

\vskip6pt
It is convenient to summarize Hayden's original proof.

\vskip6pt\noi
{\it A Sketch of Proof of Theorem 5.1}\/ \cite{Hay}.

\vskip6pt\noi
{\it Part 1}\/$.$
Let $\M$ and $\N$ be subspaces of Banach spaces $\X$ and $\Y$ such that
for every ${\phi\in\bb[\M\x\N,\FF]}$ there is ${\what\phi\in\bb[\X\x\Y,\FF]}$
for which $\what\phi|_{\M\x\N}\!=\phi$ and ${\|\what\phi\|=\|\phi\|}.$ By
Proposition 4.2, ${\bb[\M\x\N,\FF]\cong\B[\M,\N^*]}$ and
${\bb[\X\x\Y,\FF]\cong\B[\X,\Y^*]}.$ This ensures the existence of
${\what T\kern-1pt\in\B[\X,\Y^*]}$ for every ${T\kern-1pt\in\B[\M,\N^*]}$ with
${\what T|_\M\!=\kern-1ptT}$ and ${\|\what T\|\kern-1pt=\kern-1pt\|T\|}$.

\vskip6pt\noi
{\it Part 2}\/$.$
Suppose Part 1 holds for $\M=\N^*\!.$ If $T=I$, then take $\what T=E$, the
contin\-uous projection with $\R(E)=\M$ and $\|E\|=1$.

\vskip6pt\noi
{\it Part 3}\/$.$
If there is a continuous projection ${E\!:\X\to\X}$ with $\R(E)=\M$ and
$\|E\|=1$, and if ${\phi\in\bb[\M\x\Y,\FF]}$ for some Banach space $\Y$,
then $\what\phi(x,y)=\phi(Ex,y)$ for every ${(x,y)\in\X\x\Y}$ defines
$\what\phi$ in ${\bb[\X\x\Y,\FF]}$ such that
$\what\phi|_{\M\x\N}=\phi$ and $\|\what\phi\|=\|\phi\|$.

\vskip6pt\noi
{\it Part 4}\/$.$
It can be verified that Parts 2 and 3 ensure the following statement$.$ If
$\M$ is a subspace of a Banach space $\X$, and $\M=\N^*\!$ for some Banach
space $\N$, and if for every ${\phi\in\bb[\M\x\N,\FF]}$ there exists
${\what\phi\in\bb[\X\x\N,\FF]}$ such that $\what\phi|_{\M\x\N}=\phi$ and
$\|\what\phi\|=\|\phi\|$, then for every Banach space $\Y$ and every
${\phi\in\bb[\M\x\Y,\FF]}$ there exists ${\what\phi\in\bb[\X\x\Y,\FF]}$ such
that $\what\phi|_{\M\x\Y}=\phi$ and $\|\what\phi\|=\|\phi\|$,

\vskip6pt\noi
{\it Part 5}\/$.$
The statement of Theorem 5.1 can be shown to be a corollary of
Part 4$.\!\!\!$\QED

\vskip6pt
Theorem 5.2 below extends the result from \cite{Hay} restated in Theorem 5.1
by offering a necessary and sufficient condition, and showing that the
norm-one condition is required for the norm inequality only (not for the
extension) and, moreover, Theorem 5.2 holds for bilinear maps in general
(rather than for bilinear functionals).

\Theorem{5.2}
{\it Let\/ $\M$ and\/ $\N$ be linear manifolds of Banach spaces\/ $\X$ and}\/
$\Y$.

\vskip6pt\noi
{\rm(a)}
{\it Every bounded bilinear map\/ ${\phi\!:\M\x\N\to\Z}$ into an arbitrary
Banach space\/ $\Z$ has a bounded bilinear extension\/
${\what \phi\!:\X\x\Y\to\Z}$ over\/ ${\X\x\Y}$ if and only if the closures\/
$\M^-\!$ and\/ $\N^-\!$ of\/ $\M$ and\/ $\N$ are complemented subspaces of\/
$\X$ and\/ $\Y$}\/.

\vskip6pt\noi
{\rm(b)}
{\it Moreover, if\/ ${\M^-\!=\R(E)}$ and\/ ${\N^-\!=\R(P)}$ for some
projections\/ ${E\in\B[\X,\X]}$ and\/ ${P\in\B[\Y,\Y]}$ with\/
${\|E\|=\|P\|=1}$, then}\/ ${\|\what\phi\|=\|\phi\|}$.

\proof
(a) Let $\M$ and $\N$ be nonzero linear manifolds of Banach spaces $\X$ and
$\Y.$ Let $\M^-\!$ and $\N^-\!$ be their closures$.$ Take an arbitrary Banach
space $\Z$ and an arbitrary ${\phi\in\bb[\M\x\N,\Z]}.$ The following assertions
are equivalent.

\vskip6pt\noi
(a$_1$) There exists ${\what\phi\in\bb[\X\x\Y,\Z]}$ such that
${\what\phi|_{\M\x\N}\!=\phi}$.

\vskip6pt\noi
(a$_2$) $\M^-\!$ and $\N^-\!$ are complemented subspaces of $\X$ and $\Y$.

\vskip6pt\noi
(a$_1$)$\,\Rightarrow$(a$_2).$
Consider the linear manifolds $\M$ and $\N.$ Take an arbitrary nonzero bounded
linear functional ${g\!:\N\!\to\FF}$ and an arbitrary nonzero bounded linear
transformation ${T\!:\!\M\to\Z}$ (whose existences are ensured by the
Hahn--Banach Theorem), and consider the map ${\phi\!:\M\x\N\!\to\Z}$
defined by
$$
\phi(u,v)=g(v)Tu
\quad\;\hbox{for every}\;\quad
(u,v)\in\M\x\N.
$$
Since $g$ and $T$ are both linear and bounded, it is readily verified that
$\phi$ is bilinear and bounded (with $\|\phi(u,v)\|\le\|g\|\|T\|\|v\|\|u\|$
for every ${u\in\M}$ and ${v\in\N}$)$.$ Thus ${\phi\in\bb[\M\x\N,\Z]}.$ If
(a$_1$) holds, then consider its extension ${\what\phi\in\bb[\X\x\Y,\Z]}$,
and let ${\phi^v\in\B[\M,\Z]}$ and ${\what\phi^v\in\B[\X,\Z]}$ be their
$v$-sections for each ${v\in\N\sse\Y}.$ Thus for every ${u\in\M}$ and each
${v\in\N}$,
$$
\what\phi^v|_\M(u)=\what\phi(u,v)=\what\phi|_{\M\x\N}(u,v)
=\phi(u,v)=\phi^v(u).
$$
Fix an arbitrary ${v\in\N}$ for which ${g(v)\ne0}$ and
set $0\ne\alpha=g(v)\in\FF.$ Hence
$$
\what\phi^v|_\M(\cdot)=\phi^v(\cdot)=\alpha T(\cdot)\in\B[\M,\Z].
$$
Thus ${\alpha T=\phi^v\in\B[\M,\Z]}$ has an extension
${\what{\alpha T}=\what\phi^v\in\B[\X,\Z]}.$ Then there exists
${\what T\in\B[\X,\Z]}$ for which ${\what T|_\M=T}.$ Since $\X$ is a Banach
space and since this holds for every nonzero ${T\in\B[\M,\Z]}$ and every
Banach space $\Z$, then $\M^-\!$ is complemented (see, e.g.,
\cite[Theorem 3.2.17]{Meg})$.$ Similarly, by taking an arbitrary nonzero
bounded linear functional ${f\!:\!\M\to\FF}$ and an arbitrary nonzero
bounded linear transformation ${S\!:\N\!\to\Z}$, and defining
${\psi\in\bb[\M\x\N,\Z]}$ by
$$
\psi(u,v)=f(u)Sv
\quad\;\hbox{for every}\;\quad
(u,v)\in\M\x\N,
$$
the same argument ensures that $\N^-\!$ is complemented$.$ Hence (a$_2$)
holds.

\vskip6pt\noi
(a$_2$)$\,\Rightarrow$(a$_1).$
Consider the normed spaces $\M$ and $\N.$ The normed spaces ${\bb[\M\x\N,\Z]}$
and ${\B[\M,\B[\N,\Z]\kern1pt]}$ are isometrically isomorphic for every normed
space $\Z$,
$$
\bb[\M\x\N,\Z]\cong\B[\M,\B[\N,\Z]\kern1pt],
$$
by Proposition 4.2, where the natural isometric isomorphism
$$
\Im\!:\bb[\M\x\N,\Z]\to\B[\M,\B[\N,\Z]\kern1pt]
$$
that sends an arbitrary ${\phi\in\bb[\M\x\N,\Z]}$ to $T=\Im(\phi)$ in
${\B[\M,\B[\N,\Z]\kern1pt]}$ is given by $T(u)={\phi_u\in\B[\N,\Z]}$ for
each ${u\in\M}$. Therefore
$$
\phi(u,v)=T(u)(v)
\quad\;\hbox{for every}\;\quad
(u,v)\in\M\x\N.
$$
Suppose $\Z$ is a Banach space (which implies that $\B[\N,\Z]$ is a Banach
space) in order to allow extension by continuity of uniformly continuous
functions on dense sets$.$ Consider the extension by continuity
${\wtil T\in\B[\M^-\!,\B[\N,\Z]\kern1pt]}$ of
${T\in\B[\M,\B[\N,\Z]\kern1pt]}.$ For each vector ${\wtil u\in\M^-\!}$ consider
the extension by continuity ${\wtil T^{_\sim}(\wtil u)\in\B[\N^-\!,\Z]}$ of
${\wtil T(\wtil u)\in\B[\N,\Z]}$, defining a transformation
${\wtil T^{_\sim}\in\B[\M^-\!,\B[\N^-\!,\Z]\kern1pt]}$ such that
$$
\wtil T^{_\sim}(u)(v)=T(u)(v)
\quad\;\hbox{for every}\;\quad
(u,v)\in\M\x\N.
$$
Applying the above isometric isomorphic argument (i.e., Proposition 4.2) again,
$$
\bb[\M^-\!\x\N^-\!,\Z]\cong\B[\M^-\!,\B[\N^-\!,\Z]\kern1pt].
$$
Let ${\wtil\Im\!:\bb[\M^-\!\x\N^-\!,\Z]\to\B[\M^-\!,\B[\N^-\!,\Z]\kern1pt]}$
be the natural isometric isomorphism$.$ Set
${\wtil\phi=\wtil\Im^{-1}(\wtil T^{_\sim})}$ in ${\bb[\M^-\!\x\N^-\!,\Z]}$
for each $\wtil T^{_\sim}$ in ${\B[\M^-\!,\B[\N^-\!,\Z]\kern1pt]}$
so that $\wtil\phi_{\wtil u}={\wtil T^{_\sim}(\wtil u)\in\B[\N^-\!,\Z]}$
for each ${\wtil u\in\M^-\!}$, and so
$$
\wtil\phi(\wtil u,\wtil v)=\wtil T^{_\sim}(\wtil u)(\wtil v)
\quad\;\hbox{for every}\;\quad
(\wtil u,\wtil v)\in\M^-\!\x\N^-\!.
$$
Thus $\wtil\phi(u,v)=\wtil T^{_\sim}(u)(v)=T(u)(v)=\phi(u,v)$
for every ${(u,v)\kern-1pt\in\!\M\x\N}.$ Therefore
$$
\wtil\phi|_{\M\x\N}=\phi.
$$
Now suppose (a$_2$) holds$.$ By Proposition 4.1 this means there are
projections $E$ in ${\B[\X,\X]}$ and $P$ in ${\B[\Y,\Y]}$ with
${\R(E)=\M^-\!}$ and ${\R(E)=\N^-\!}.$ Set
$$
\what\phi(x,y)=\wtil\phi(Ex,Py)
\quad\;\hbox{for every}\;\quad
(x,y)\in\X\x\Y.
$$
As ${E\!:\!\X\!\to\X}$ and ${P\!:\!\Y\to\!\Y}$ are linear and bounded, and
${\wtil\phi\!:\!\M^-\!\x\N^-\!\to\kern-1pt\Z}$ is bilinear and bounded, then
${\what\phi\!:\X\x\Y\to\kern-1pt\Z}$ is bilinear and bounded (i.e.,
${\what\phi\in\bb[\X\x\Y,\Z]}$)$.$ Also, as $E$ and $P$ act as
the identity on ${\M\sse\R(E)}$ and ${\N\sse\R(P)}$, we get (a$_1$)$:$
$$
\what\phi|_{\M\x\N}=\wtil\phi(E|_\M,P|_\N)=\wtil\phi|_{\M\x\N}=\phi.
$$

\vskip2pt\noi
(b)
Hence
$\|\phi\|=\|\what\phi|_{\M\x\N}\|
\le\|\what\phi\|
=\|\wtil\phi(E(\cdot),P(\cdot))\|
\le\|\wtil\phi\|\kern.5pt\|E\|\kern.5pt\|P\|.$
Moreover, although there is no extension by continuity for bilinear maps,
\begin{eqnarray*}
\|\wtil\phi\|
&\kern-6pt=\kern-6pt&
\!\!\sup_{\!0\ne\wtil u\in\M^-\!,\,0\ne\wtil v\in\N^-\!\!}\!
\smallfrac{\|\phi(\wtil u,\wtil v)\|}{\|\wtil u\|\|\wtil v\|}
=\!\!\sup_{\!0\ne\wtil u\in\M^-\!,\,0\ne\wtil v\in\N^-\!\!}\!\!\!
\smallfrac{\|\wtil T^{_\sim}(\wtil u)(\wtil v)\|}{\|\wtil u\|\|\wtil v\|}
=\!\!\sup_{\!0\ne u\in\M,\,0\ne v\in\N\!}\!\!\!
\smallfrac{\|\wtil T^{_\sim}(u)(v)\|}{\|u\|\|v\|}
\\
&\kern-6pt=\kern-6pt&
\sup_{\!0\ne u\in\M,\,0\ne v\in\N\!}
\smallfrac{\|T(u)(v)\|}{\|u\|\|v\|}
=\sup_{\!0\ne u\in\M,\,0\ne v\in\N\!}
\smallfrac{\|\phi(u,v)\|}{\|u\|\|v\|}
=\|\phi\|.
\end{eqnarray*}
Thus $\|\phi\|\le\|\what\phi\|\le\|\phi\|\kern.5pt\|E\|\kern.5pt\|P\|.$ Then
$\|E\|=\|P\|=1$ implies $\|\what\phi\|=\|\phi\|$. \QED

\vskip6pt
$\!\!$The converse of Theorem 5.2(b) fails:
$\!{\|\what\phi\|\kern-1pt=\kern-1pt\|\phi\|}\kern-.5pt$ does not imply
$\kern-.5pt{\|E\|\kern-1pt=\kern-1pt\|P\|\kern-1pt=\kern-1pt1\kern-.5pt}.$ For
instance, ${E\in\B[\RR^2\!,\RR^2]}$ given by ${E(x_1,x_2)}={(0,x_1+x_2)}$
defines a projection on $\RR^2$ with ${\|E\|=\sqrt2}.$ Take
${f\in\B[\RR^2,\RR]}$ given by ${f((x_1,x_2))=x_2}$ for ${(x_1,x_2)\in\RR^2}$
so that $\|f\|=1.$ Take ${\what\phi\in\bb[\RR^2\x\RR^2\!,\RR]}$ given by
${\what\phi(x,y)=f(x)\kern.5ptf(y)}$ for ${(x,y)\in\RR^2\x\RR^2}\!$ so that
$\|\what\phi\|=\|f\|^2.$ Set ${\M=\R(E}).$ The restriction
${\phi=\what\phi|_{\M\x\M}\in\bb[\M\x\M,\RR]}$ is such that
$\phi((0,\alpha),(0,\beta))=\alpha\beta$ for
${(0,\alpha),\,(0,\beta)\in\kern-1pt\M}$ and so ${\|\phi\|=1}.$

\vskip6pt
If $\X$ and $\Y$ are Hilbert spaces, then one gets a full extension with no
restrictions on the linear manifolds $\M$ and $\N$, as expected.

\Corollary{5.3}
{\it Every bounded bilinear map\/ ${\phi\!:\!\M\x\N\!\to\kern-1pt\Z}$ defined
on the \hbox{Cartesian} product of arbitrary linear manifolds\/ $\M$ and\/
$\N$ of arbitrary Hilbert spaces\/ $\X$ and\/ $\Y$ into an arbitrary Banach
space\/ $\Z$ has a bounded bilinear extension\/ ${\what\phi\!:\X\x\Y\to\Z}$
over\/ ${\X\x\Y}$ such that}\/ ${\|\what\phi\|=\|\phi\|}$.

\proof
(a) Hilbert spaces are complemented$.$ Thus every subspace of a Hilbert
space is complemented, and so Theorem 5.2(a) applies to every linear manifold
$\M$ of a Hilbert space $\X$ and to every linear manifold $\N$ of a Hilbert
space $\Y$.

\vskip6pt\noi
(b) If $\X$ ad $\Y$ are Hilbert spaces, then the orthogonal projections
${E\in\B[\X,\X]}$ with ${\R(E)=\M^-\!}$ and ${P\in\B[\Y,\Y]}$ with
${\R(P)=\N^-\!}$ are such that $\|E\|=\|P\|=1$, and therefore
${\|\what\phi\|=\|\phi\|}$ by Theorem 5.2(b). \QED

\vskip6pt
Extensions of linear transformations (or forms) are not unique, and so
extensions of bilinear maps are not unique (since a product of linear forms
is a bilinear form).

\Remark{5.4}
Consider the following classes of operators on a Banach space $\X$.
$$
\Gamma_R[\X]
=\big\{T\in\B[\X,\X]\!:\,\R(T)^-\;
\hbox{is a complemented subspace of $\X$}\big\},
$$
\vskip-4pt\noi
$$
\Gamma_N[\X]
=\big\{T\in\B[\X,\X]\!:\,\N(T)\;
\hbox{is a complemented subspace of $\X$}\big\},
$$
$$
\Gamma[\X]=\Gamma_R[\X]\cap\Gamma_N[\X],
\qquad
\Se\Gamma[\X]=\Gamma_R[\X]\cup\Gamma_N[\X].
$$
The class $\Gamma[\X]$ is large enough$.$ It includes, for instance, the class
of all compact perturbations of semi-Fredholm operators (see, e.g.,
\cite[Remark 5.1(b)]{Kub1})$.$ A straightforward consequence of Theorem 5.2
leads to the following result.
\vskip6pt\noi
{\narrower
{\it Let\/ $\M$ and\/ $\N$ be linear manifolds of\/ Banach spaces $\X$ and\/
$\Y$ for which there are operators ${T\in\Se\Gamma[\X]\sse\B[\X,\X]}$ and\/
${S\in\Se\Gamma[\Y]\sse\B[\Y,\Y]}$ such that
$$
\M^-\!=\R(T)^-\!
\;\;\hbox{or}\;\;
\M^-\!=\N(T)
\;\quad\;\hbox{and}\;\quad
\N^-\!=\R(S)^-\!
\;\;\hbox{or}\;\;
\N^-\!=\N(S),
$$
according to whether\/ $T$ lies in\/ $\Gamma_R[\X]$ or\/ $\Gamma_N[\X]$ and\/
$S$ lies in\/ $\Gamma_R[\Y]$ or\/ $\Gamma_N[\Y].$ If\/
${\phi\!:\M\x\N\!\to\kern-1pt\Z}$
is a bounded bilinear map into an arbitrary Banach space\/ $\Z$, then there
exists a bounded bilinear extension\/ ${\what\phi\!:\X\x\Y\to\kern-1pt\Z}$
of\/ $\phi$ over}\/ ${\X\x\Y}$.
\vskip2pt}

\section{An Application to Projective Tensor Products}

If $\M$ and $\N$ are linear manifolds of linear spaces $\X$ and $\Y$, then
${\M\otimes\N}$ is a regular linear manifold of the linear space
${\X\otimes\Y}.$ If ${\X\otimes\Y}$ is equipped with the injective norm, then
${\M\otimesv\N}$ is a linear manifold of ${\X\otimesv\Y}.$ However, this is
not always the case if ${\X\otimes\Y}$ is equipped with the projective norm$.$
The next corollaries give necessary and sufficient conditions for
${\M\otimesw\N}$ to be a linear manifold of ${\X\otimesw\Y}.$ Moreover,
extensions of bounded linear transformations on regular linear manifolds of
projective tensor products also comes as another consequence of Theorem 5.2.

\Corollary{6.1}
{\it Let\/ ${\X,\Y,\Z}$ be arbitrary Banach spaces, let\/ $\M$ and\/ $\N$ be
subspaces of\/ $\X$ and\/ $\Y$, respectively, and consider the following
assertions}\/.
\begin{description}
\item{$\kern-12pt$\rm(o)}
{\it $\;\M$ and\/ $\N$ are complemented in\/ $\X$ and\/ $\Y$ with\/ $\M=\R(E)$
and\/ $\N=\R(P)$ for projections\/ ${E\in\B[\X,\X]}$ and\/ ${P\in\B[\Y,\Y]}$
such that}\/ ${\|E\|=\|P\|=1}$.
\vskip6pt
\item{$\kern-12pt$\rm(a)}
{\it $\;{\M\otimesw\N}$ is a linear manifold of}\/ ${\X\otimesw\Y}$.
\vskip6pt
\item{$\kern-12pt$\rm(b)}
{\it $\,$Every\/ bounded bilinear map\/ ${\phi\!:\M\x\N\!\to\Z}$ has a bounded
bilinear extension\/ ${\what\phi\!:\X\x\Y\to\Z}$ with}\/
${\|\what\phi\|=\|\phi\|}$.
\vskip6pt
\item{$\kern-12pt$\rm(c)}
{\it $\;$Every bounded linear transformation\/ ${T\!:\M\otimesw\N\!\to\Z}$
has a \hbox{$\;$bounded} linear extension\/ ${\wtil T\!:\X\otimesw\Y\to\Z}$
with}\/ ${\|\wtil T\|=\|T\|}$.
\vskip6pt
\item{$\kern-12pt$\rm(d)}
{\it $\,$Every bounded linear transformation\/
${\overline T\!:\M\hotimesw\N\to\Z}$ has a \hbox{$\;$bounded} linear
extension\/ ${\what T\!:\X\hotimesw\Y\to\Z}$ with}\/
${\|\what T\|=\|\overline T\|}$.
\end{description}
{\it The above assertions are related as follows}\/:
$$
\hbox{{\rm(o)} $\limply$ {\rm(a,b)}},
\qquad
\hbox{{\rm(a,b)}$\iff${\rm(c)}$\iff${\rm(d)}}.
$$

\proof
Let\/ ${\X,\Y,\Z}$ be arbitrary Banach spaces.

\vskip6pt\noi
{\rm(o)}$\;\limply\,${\rm(a,b)}$.$
If (o) holds, then according to Proposition 4.3 and Theorem 5.2,
\vskip4pt\noi
\begin{description}
\item{\rm(a)}
$\;{\M\otimesw\N}$ is a linear manifold of ${\X\otimesw\Y}$
\vskip0pt
(or, equivalently, ${\M\otimesw\N}$ is a subspace of ${\X\otimesw\Y}$),
\quad and
\vskip4pt
\item{\rm(b)}
$\;$every ${\phi\in\bb[\M\x\N,\Z]}$ into an arbitrary Banach space
$\Z$ has an extension ${\quad\,\what\phi\in\bb[\X\x\Y,\Z]}$ with
${\|\what\phi\|=\|\phi\|}$.
\end{description}
\vskip4pt\noi
{\rm(a,b)}$\iff${\rm(c)}$\iff${\rm(d)}$.$
Each extension in (c) or (d) only makes sense if (a) holds$.$ According to
Proposition 4.4
$$
\B[\X\otimesw\Y,\Z]\cong\bb[\X\x\Y,\Z].
$$
The natural isometric isomorphism between them,
$$
\J\!:\B[\X\otimesw\Y,\Z]\to\bb[\X\x\Y,\Z]
$$
such that ${\psi=\J\kern.5pt(S)}$ for ${S\in\B[\X\otimesw\Y,\Z]}$ and
${\psi\in\bb[\X\x\Y.\Z]}$, is given by
$$
\psi(x,y)=\J\kern.5pt(S)(x,y)=S(x\otimes y)
\quad\,\hbox{and so}\,\quad
S(x\otimes y)=\J^{-1}(\psi)(x\otimes y)=\psi(x,y)
$$
for every ${(x,y)\in\X\x\Y}.$ By Theorem 5.2(a) $\M$ and $\N$ are
complemented if and only if every $\phi$ in ${\bb[\M\x\N,\Z]}$ has an
extension $\what\phi$ in ${\bb[\X\x\Y,\Z]}.$ Thus in this case with $T$ in
${\B[\M\otimesw\N,\Z]}$ and $\wtil T$ in ${\B[\X\otimesw\Y,\Z]}$ being the
isometrically isomorphic images of $\phi$ in ${\bb[\M\x\N,\Z]}$ and
$\what\phi$ in ${\bb[\X\x\Y,\Z]}$,
$$
T\cong\phi=\what\phi|_{\M\x\N}
\quad\;\hbox{and}\;\quad
\what\phi\cong\wtil T.
$$
As the restriction ${\wtil T|_{\M\otimesw\N}\in\B[\M\otimesw\N,\Z]}$
of ${\wtil T\in\B[\X\otimesw\Y,\Z]}$ to ${\M\otimesw\N}$ only makes sense if
(a) holds, then in this case with ${\phi=\J\kern.5pt(T)}$ and
${\what\phi=\J\kern.5pt(\wtil T)}$,
\begin{eqnarray*}
T(u\otimes v)
&\kern-6pt=\kern-6pt&
\J^{-1}(\phi)(u\otimes v)
=\phi(u,v)=\what\phi|_{\M\x\N}(u,v)=\what\phi(u,v)
=\J(\wtil T)(u,v)                                                        \\
&\kern-6pt=\kern-6pt&
\wtil T(u\otimes v)
=\wtil T|_{\M\otimesw\N}(u\otimes v)
\quad\;\hbox{for every}\;\quad
(u\otimes v)\in\M\otimesw\N.
\end{eqnarray*}
Therefore ${T(\digamma)=\wtil T|_{\M\otimesw\N}(\digamma)}$ for every
${\digamma\kern-1pt=\kern-1pt\sum_iu_i\otimes v_i\in\M\otimesw\N}.$ Hence
$$
\what\phi|_{\M\x\N}=\phi
\quad\;\hbox{implies}\;\quad
\wtil T|_{\M\otimesw\N}=T.
$$
\vskip0pt\noi
{\narrower
Thus
{\it if every\/ ${\phi\in\bb[\M\x\N,\Z]}$ has an extension\/
${\what\phi\in\bb[\X\x\Y,\Z]}$ and\/ {\rm(a)} holds, then every\/
${T\in\B[\M\otimesw\N,\Z]}$ has an extension}\/
${\wtil T\in\B[\X\otimesw\Y,\Z]}$.
\vskip6pt}
\vskip0pt\noi
A symmetric argument ensures the converse:
\vskip6pt\noi
{\narrower
{\it if\/ every\/ ${T\in\B[\M\otimesw\N,\Z]}\,$ has an extension\/
${\wtil T\in\B[\X\otimesw\Y,\Z]}$, then every\/ ${\phi\in\bb[\M\x\N,\Z]}$
has an extension\/ ${\what\phi\in\bb[\X\x\Y,\Z]}$ and {\rm(a)} holds}\/.
\vskip6pt}
\vskip0pt\noi
Suppose either (c) or (d) holds, and so (a) holds$.$ Thus the completions
${\M\hotimesw\N}$ of ${\M\otimesw\N}$ and ${\X\hotimesw\Y}$ of
${\X\otimesw\Y}$ are such that ${\M\hotimesw\N}$ is a subspace of
${\X\hotimesw\Y}.$ Since $\Z$ is a Banach space,
$$
\B[\M\hotimesw\N,\Z]\cong\B[\M\otimesw\N,\Z]
\quad\;\hbox{and}\;\quad
\B[\X\hotimesw\Y,\Z]\cong\B[\X\otimesw\Y,\Z].
$$
Again, let $\overline T$ in ${\B[\M\hotimesw\N,\Z]}$ be the extension over
completion of $T$ in ${\B[\M\otimesw\N,\Z]}$ so that ${\overline T\cong T}$,
and let $\what T$ in ${\B[\X\hotimesw\Y,\Z]}$ be the extension over completion
of $\wtil T$ in ${\B[\X\otimesw\Y,\Z]}$ so that ${\what T\cong\wtil T}.$
Extensions over completions are unique up to isometric isomorphism$.$ Then
$\what T$ in ${\B[\X\hotimesw\Y,\Z]}$ extends $\overline T$ in
${\B[\M\hotimesw\N,\Z]}$ if and only if $\wtil T$ in ${\B[\X\otimesw\Y,\Z]}$
extends $T$ in ${\B[\M\otimesw\N,\Z]}$:
$$
\wtil T|_{\M\otimesw\N}\cong\what T|_{\M\hotimesw\N}=\overline T\cong T
\quad\;\hbox{and}\;\quad
\what T|_{\M\hotimesw\N}\cong\wtil T|_{\M\otimesw\N}=T\cong \overline T.
$$
\vskip0pt\noi
{\narrower
Hence
{\it every\/ ${\overline T\kern-1pt\in\kern-1pt\B[\M\hotimesw\N,\Z]}$ has an
extension\/ ${\what T\kern-1pt\in\kern-1pt\B[\X\hotimesw\Y,\Z]}$ if and only
if every\/ ${T\in\B[\M\otimesw\N,\Z]}$ has an extension}\/
${\wtil T\in\B[\X\otimesw\Y,\Z]}$.
\vskip6pt}
\vskip0pt\noi
Also, since ${\what\phi\cong\wtil T}$, $\,{\phi\cong T}$,
$\,{\wtil T\cong\what T}$, and ${T\cong \overline T}$, if one of the norms
(b), (c), or (d) coincide, then so does the others. \QED

\vskip6pt
A particular case with a rather simplified statement is immediately obtained
by fixing ${\Z=\FF}$ in Corollary 6.1 as follows$.$ (This extends
\cite[Proposition 2.11]{Rya}).

\Corollary{6.2}
{\it If\/ $\M$ and $\N$ are subspaces of Banach spaces $\X$ and $\Y$, then
the following assertions are pairwise equivalent}\/.
\begin{description}
\item{$\kern-12pt$\rm(a)}
{\it $\,{\M\otimesw\N}$ is a linear manifold of}\/ ${\X\otimesw\Y}$.
\vskip6pt
\item{$\kern-12pt$\rm(b)}
{\it $\,$Every\/ bounded bilinear form\/ ${\phi\!:\M\x\N\!\to\FF}$
has a bounded bilinear extension\/ ${\what\phi\!:\X\x\Y\to\FF}$ with}\/
${\|\what\phi\|=\|\phi\|}$.
\vskip6pt
\item{$\kern-12pt$\rm(c)}
{\it $\,$Every bounded linear form\/ ${f\!:\M\otimesw\N\to\FF}$ has a
bounded linear extension\/ ${\wtil f\in\X\otimesw\Y}$ with}\/
${\|\wtil f\|=\|f\|}$.
\vskip6pt
\item{$\kern-12pt$\rm(d)}
{\it $\,$Every bounded linear form \/ ${\overline f\!:\M\hotimesw\N\to\FF}$
has a bounded linear extension\/ ${\what f\!:\X\hotimesw\Y\to\FF}$
with}\/ ${\|\what f\|=\|\overline f\|}$.
\end{description}

\proof
If (b) holds, then the set of all bilinear forms in ${\bb[\M\x\N,\FF]}$ with
norm less than 1 is included in the set of all restrictions to ${\M\x\N}$ of
all bilinear forms in ${\bb[\X\x\Y,\FF]}$ with norm less that 1$.$ Then by
Proposition 4.5
$$
\|\digamma\|_{\M\otimesw\N}
\le
\sup_{\|\psi\|\le1,\,\psi\in\sbb[\X\x\Y,\FF]}\Big|{\sum}_i\psi(u_i,v_i)\Big|
=\|\digamma\|_{\X\otimesw\Y}.
$$
Since ${\|\digamma\|_{\X\otimesw\Y}\le\|\digamma\|_{\M\otimesw\N}}$, then
${\|\digamma\|_{\M\otimesw\N}=\|\digamma\|_{\X\otimesw\Y}}$ which is
equivalent to (a)$.$ So (b) implies (a)$.$ Conversely, (a) implies (c) by the
Hahn--Banach Theorem$.$ Also, (c) implies (b) and (c) is equivalent to (d) by
Corollary 6.1. \QED

\vskip6pt
According to Corollary 5.3, in a Hilbert-space setting each assertion in
Corollaries 6.1 and 6.2 holds true.

\bibliographystyle{amsplain}

\end{document}